\documentclass{ifacconf}
\usepackage{natbib}        
\usepackage{amssymb, amsmath,amsfonts,verbatim, mathtools,enumerate,algorithm,dsfont}
\usepackage{comment}
\usepackage{diagbox}
\usepackage{times,algpseudocode}

\allowdisplaybreaks[2]
\newcommand{\ones}{\mathbf 1}
\newcommand{\reals}{{\mathbb{R}}}
\newcommand{\naturals}{{\mathbb{N}}}

\newcommand{\argmin}{\mathop{\rm argmin}}

\newcommand{\norm}[1]{\left\lVert#1\right\rVert}
\newcommand{\mnorm}[1]{{\left\vert\kern-0.25ex\left\vert\kern-0.25ex\left\vert #1 
    \right\vert\kern-0.25ex\right\vert\kern-0.25ex\right\vert}}

\newcommand{\mc}{\mathcal}

\algnewcommand\algorithmicforeach{\textbf{for each}}
\algdef{S}[FOR]{ForEach}[1]{\algorithmicforeach\ #1\ \algorithmicdo}

\newcommand{\cl}{\operatorname{cl}}
\newcommand{\diam}[1]{\operatorname{diam}\left(#1\right)}
\usepackage{graphicx}      

\usepackage[dvipsnames]{xcolor}

 \newcommand\ForAuthorsAssa[1]
 {\par\smallskip                     
  \begin{center}
  \fbox
  {\parbox{0.9\linewidth}
    {\raggedright\sc--- {\color{blue}#1}}
  }
  \end{center}
  \par\smallskip                     %
 }
\begin{document}
\begin{frontmatter}

\title{Fixed Points of Set-Based Bellman Operator} 


\author[SL]{Sarah H.Q. Li} 
\author[AA]{Assal{\'e}  Adj{\'e}} 
\author[Third]{Pierre-Lo{\"\i}c Garoche}
\author[SL]{Beh\c cet A\c c\i kme\c se}

\address[SL]{William E.Boeing Department of Aeronautics and Astronautics, University of Washington, Seattle, USA. (e-mail: sarahli@uw.edu, behcet@uw.edu).}
\address[AA]{LAMPS, Universit{\'e} de Perpignan Via Domitia, Perpignan, France
 (e-mail: assale.adje@univ-perp.fr)}
\address[Third]{ONERA -- The French Aerospace Lab, Univ. of Toulouse, France, (e-mail: pierre-loic.garoche@onera.fr)}
\begin{abstract}                
Motivated by uncertain parameters encountered in Markov decision processes (MDPs), we study the effect of parameter uncertainty on Bellman operator-based methods. 
Specifically, we consider a family of MDPs where the cost parameters are from a given compact set. We then define a Bellman operator acting on an input set of value functions to produce a new set of value functions as the output under all possible variations in the cost parameters. Finally we prove the existence of a {\em fixed point} of this set-based Bellman operator by showing that it is a contractive operator on a complete metric space.  
\end{abstract}

\end{frontmatter}
\section{Introduction}

Markov decision process (MDP) is a widely used mathematical framework for control design in stochastic environments, eg. density control of a swarm of agents~\citep{accikmecse2012markov, demir2015decentralized}. It is also a fundamental framework for reinforcement learning, robotic motion planning and stochastic games~\citep{filar2012competitive,li2019tolling}. An MDP can be solved for different objectives including  minimum average cost, minimum discounted cost, reachability, among others~\citep{puterman2014markov}. Given an objective, solving an MDP is equivalent to computing the optimal policy and  the optimal \emph{value function} of a decision maker over the state space. Among different algorithms for computing the optimal policy, most are based on the Bellman equation that characterizes the value function as its fixed point. 

In applications, it is common to encounter MDPs with uncertainties. When modeling an environment as a stochastic process, sampling techniques are often used to determine process parameters such as process costs or probabilities; such models are inherently uncertain. In stochastic games, the cost and probability parameters change with respect to another decision maker's strategy. While existing works focus on certain perturbations in MDPs~\citep{bielecki1991singularly,altman1993stability, abbad1992perturbation},  these results  do not generalize to the analysis of overall behaviour of the MDP under all possible cost parameters in a compact set. 

Additionally, how uncertainty in MDP cost parameters affect the outcome of value iteration type methods is not well studied. Dynamic programming on bounded MDPs is studied in~\cite{givan2000bounded} for specifically interval sets, however convergence over general compact sets is not considered.  While computation of the fixed points of Bellman operator is the topic of numerous studies~\citep{delage2010percentile}, most focus on the convergence analysis of value iteration  and its stopping criteria~\citep{ashok2017value, eisentraut2019stopping}. However, they do not consider the relationship between bounds on the optimal value function and the uncertainty in cost. Similarly motivated, \cite{haddad2018interval} analyzes entry-wise uncertain transition kernels by using graph-based MDP transformations. While we also derive bounds of an MDP due to uncertain parameters, we differ in our approach: our set-based framework allows for direct extraction of the value iteration trajectories with respect to the set of cost parameters. This differentiates our work from~\cite{haddad2018interval} due to their graphical abstraction of MDP, which allows for derivation of bounds but not extraction of value function trajectories. 

\textbf{Contributions:} We characterize the solutions of a family of MDPs at once, represented as \emph{sets} of MDPs. More specifically, we: (i) develop a characterization of MDPs with uncertain cost parameters; (ii) propose a set-based Bellman operator over non-empty compact sets; (iii)  establish  the contractivity of this set-based Bellman operator  with the existence of a unique compact fixed point set.


\section{Review of MDPs and Bellman Operator}
\label{sec:setup}

\textbf{Notation}: Sets of $N$ elements are given by $[N] = \{0, \ldots, N-1\}$. We denote the set of matrices of $i$ rows and $j$ columns with real (non-negative) valued entries as $\reals^{i \times j}(\reals_+^{i\times j})$. Elements of sets and matrices are denoted by capital letters, $X$, while sets are denoted by cursive letters, $\mc{X}$. The ones column vector is denoted by $\ones_N = [1, \ldots, 1]^T \in \reals^{N\times 1}$.

\subsection{MDP}
We consider a \emph{discounted infinite-horizon MDP} defined by $([S], [A], P, C, \gamma)$ for a decision maker, where
\begin{enumerate}
	\item $[S]$ denotes the finite set of states.

	\item $[A]$ denotes the finite set of actions. Without loss of generality, assume that every action is admissible from each state $s \in [S]$.

	\item $P \in \reals^{S \times SA}$ denotes the transition kernel. Each component $P_{s',sa}$ is the probability of arriving in state $s'$ by taking state-action $(s,a)$. Matrix $P$ is column stochastic and element-wise non-negative --- i.e. $\sum_{(s,a) \in [S]\times[A]} P_{s',sa} = 1$, $P_{s',sa} \geq 0 $, $\forall \ s',s \in [S], a \in [A]$.
	
	\item $C \in \reals^{S \times A}$ denotes the cost of each pair $(s,a)$. 
	
	\item $\gamma \in (0,1)$ denotes the discount factor. 
\end{enumerate}

At each time step $t$, the decision maker chooses an action $a$ based on its current state $s$. The state-action pair $(s, a)$ induces a probability distribution $P_{(\cdot), sa} \in \reals^S$, where $P_{s',sa}$ is the probability that the decision maker arrives at $s'$ at time step $t+1$. The state-action $(s,a)$ also induces a cost $C_{sa}$ that must be paid by the decision maker. 

At each time step, the decision maker chooses a \emph{policy} that dictates the action chosen at each state $s$. We denote policy as a function $\pi:S \times A \rightarrow \reals_+$, where $\pi(s,a)$ denotes the probability that action $a$ is chosen at state $s$. 
We denote the set of all feasible policies of an MDP by $\Pi$. In our context, it suffices to consider only \emph{deterministic, stationary} policies i.e. $\pi(s,a)$ is a time invariant function that returns $1$ for exactly one action, and $0$ for all other possible actions.

We denote the policy matrix induced by a policy as $M_\pi \in \reals^{S\times SA}$, where \[(M_\pi)_{s', sa} = \begin{cases}
\pi(s,a) & s' = s \\
0 & s' \neq s
\end{cases}.
\]
Every stationary policy induces a stationary \emph{Markov chain} \citep{el2018controlled}, given by $PM_{\pi}^T$. 
Each stationary policy also induces a stationary cost given by
\begin{equation}\label{eqn:policyCost}
  C(\pi) = \sum_{i \in [S]} e_ie_i^T M_{\pi}(\ones_S \otimes I_A) C^Te_i, \ C(\pi) \in \reals^{S},  
\end{equation}
where $e_i \in \reals^{S}$ is the unit vector pointing in the $i^{th}$ coordinate.

For an MDP $([S],[A],P,C,\gamma)$, we are interested in minimizing the \emph{discounted infinite horizon expected cost}, defined with respect to a policy $\pi$ as
\begin{equation}\label{eqn:expectationOptimalV}
 V^\star_{s^0} = \min_{\pi \in \Pi} \mathbb{E}^\pi_{s^0} \Big\{ \sum_{t = 0}^\infty \gamma^t C_{s^t a^t}\Big\}, \quad \forall \ s_0 \in [S]
\end{equation}
where $\gamma \in (0,1)$ is the discount factor of future cost, $s^t$ and $a^t$ are the state and action taken at time step $t$,  and $s^0$ is the state that the decision maker starts from at $t = 0$. 

The minimum expected cost $V_s^\star$ is called the \emph{optimal value function}. The policy $\pi^\star$ that achieves the optimal value function is called an \emph{optimal policy}. In general, $V^\star_s$ is unique while $\pi^\star$ is not. It is well known that the set of optimal policies always includes at least one \emph{deterministic} stationary policy~\citep[Thm 6.2.11]{puterman2014markov} --- i.e. for each $s$, $\pi(s,a)$ returns $1$ for exactly one action, and $0$ for all other possible actions.

\subsection{Bellman Operator}
Determining the optimal value function of a given MDP is equivalent to solving for the fixed point of the associated Bellman operator, for which a myriad of techniques exists~\citep{puterman2014markov}. We introduce the Bellman operator here as well as relating its fixed point to the corresponding MDP problem.

\begin{defn}[Standard Bellman Operator]\label{def:bellmanOp} 
For a discounted infinite horizon MDP $([S],[A],P,C,\gamma)$, its associated Bellman operator $f_C: \reals^{S}\rightarrow \reals^{S}$ is given component-wise by 
\[\Big(f_C(V)\Big)_s :=  \min_{a \in [A]}\ C_{sa} +\gamma\sum_{s' \in [S]} P_{s'sa}V_{s'},\ \forall\, s\in [S]. \]
\end{defn}
The fixed point of the Bellman operator is a value function $V \in \reals^S$ that is invariant under the operator. 
\begin{defn}[Fixed Point]\label{def:fixedPoint}
Let $F:\mc{X} \to \mc{X}$ be an operator on the metric space $\mc{X}$, $ V^\star \in \mc{X}$ is a fixed point of $F$ if it satisfies 
\begin{equation}
V^\star = F(V^\star).
\end{equation}
\end{defn}

In order to show that the Bellman operator has a unique fixed point, we consider the following operator property. 
\begin{defn}[Contraction Operator]\label{def:contractOp}
Let $(\mc{X},d)$ be a complete metric space. An operator $F:\mc{X} \to \mc{X}$ is a contraction operator if it satisfies
\[d(F(\mc{V}), F(\mc{V}')) < d(\mc{V}, \mc{V}'), \quad \forall \ \mc{V}, \ \mc{V}' \ \in \mc{X}.\]
\end{defn}
The Bellman operator is known be a contraction operator on the complete metric space $(\reals^{S}, \norm{\cdot}_\infty)$. From the Banach fixed point theorem~\citep{puterman2014markov}, it has a unique fixed point. Because the optimal value function $V^\star$ is given by the unique fixed point of the associated Bellman operator, we use the terms optimal value function and fixed point of $f_C$ interchangeably. 

In addition to obtaining $V^\star$, MDPs are also solved to determine the \emph{optimal policy}, $\pi^\star$. We note that because every feasible policy $\pi$ induces a Markov chain, $\pi$ also induces a unique stationary value function $V$ which satisfies
\begin{equation}\label{eqn:vecForm_of_VI}
    V = C(\pi) + \gamma M_{\pi}P^TV.
\end{equation}
Given a feasible policy $\pi$, we can equivalently solve for the stationary value function $V$ as $V = (I - \gamma M_{\pi}P^T )^{-1}C(\pi)$. From this perspective, the optimal value function is the minimum vector among the finite set of stationary value functions generated by the set of all policies $\Pi$. 

From the optimal value function $V^\star$, we can also derive a deterministic optimal policy from the Bellman operator as
\begin{equation}\label{eqn:optimalPol}
    \pi^\star(s, a) = \begin{cases} 1 & a = \underset{\bar{a} \in [A]}{\argmin} \ C_{s\bar{a}} + \gamma \underset{s'\in [S]}{\sum}P_{s',s\bar{a}}V^\star_{s'} \\
    0 & \text{otherwise} 
    \end{cases},\ \forall \ s \in [S].
\end{equation}
While the optimal policy does not need to be deterministic and stationary, the optimal policy $\pi^\star$ derived from~\eqref{eqn:optimalPol} will always be deterministic.

\subsection{Termination Criteria for Value Iteration}
Among different algorithms to determine the fixed point of the Bellman operator, \emph{value iteration} (VI) is a commonly used and simple technique in which the Bellman operator is iteratively applied until the optimal value is reached --- i.e. starting from any value function $V^0 \in \reals^S$ and $k = 1, \ldots$, we apply
\begin{equation} \label{eqn:VI}
    V^{k+1}_s = \min_{a \in [A]} \ C_{sa} + \gamma \sum_{s'\in [S]} P_{s',sa}V^k_{s'}, \quad \forall s \in [S].
\end{equation}
The iteration scheme given by~\eqref{eqn:VI} converges to the fixed point of the corresponding discounted infinite horizon MDP. 
The \emph{stopping criteria} of VI can be considered the over-approximation of the optimal value function.
\begin{lem}\label{thm:stoppingCriteria}\citep[Thm. 6.3.1]{puterman2014markov}
For any initial value function $V^0 \in \reals^S$, let $\{V^k\}_{k\in \naturals}$ be the value function trajectory from~\eqref{eqn:VI}. Whenever there exists $\epsilon > 0$, such that
$\norm{V^{k+1} - V^k} < \epsilon \frac{(1 - \gamma)}{2\gamma}$, then $V^{k+1}$ is within $\epsilon /2$ of the fixed point $V^\star$, i.e. 
$\norm{V^{k+1} - V^\star} < \frac{\epsilon}{2}$. 
\end{lem}
Lemma~\ref{thm:stoppingCriteria} connects the sequence $\{V^k\}_{k \in \naturals}$'s relative convergence to its absolute convergence towards $V^\star$ by showing that the former implies the latter. In general, the stopping criteria differ for different MDP objectives (see~\cite{haddad2018interval} for recent results on stopping criteria for reachability). 

\section{Set-based Bellman Operator}
\label{sec:set-VI}
The standard  Bellman operator with respect to a fixed cost parameter $C$ is well studied. Motivated by a \emph{family} of MDPs corresponding to a compact set of cost parameters  $\mc{C} \subseteq \reals^{S\times A}$ with all other data parameters remaining identical, we \emph{lift} the Bellman operator to operate on sets rather than individual vectors in $\reals^S$. For the set-based operator, we analyze its set-based domain and prove that it is a contraction operator. We also prove the existence of a unique fixed point \emph{set} $\mc{V}^\star$ for a set-based Bellman operator and relate its properties to the fixed point of the standard Bellman operator.  
\subsection{Set-based operator properties}
We define a new metric space $(H(\reals^{S}), d_H)$ based on the Banach space $(\reals^S, \norm{\cdot}_\infty)$ to serve as our set-based operator domain~\citep{rudin1964principles}, where $H(\reals^{S})$ is the collection of non-empty compact subsets of $\reals^{S}$ equipped with \emph{partial order}: for $\mc{V}, \mc{V}' \in H(\reals^S)$, $\mc{V} \preceq \mc{V}'$ if $\mc{V} \subseteq \mc{V}'$ --- i.e. if $\mc{V}$ is a subset of $\mc{V}'$. The metric $d_H$ is the \emph{Haussdorf distance}~\citep{henrikson1999completeness} defined as 
\begin{equation}
\begin{aligned}
d_H(\mc{V}, \mc{V}') = \max\{ &\sup_{V \in \mc{V}}\inf_{V' \in \mc{V}'} \norm{V - V'}_\infty,\\
&\sup_{V' \in \mc{V}'}\inf_{V \in \mc{V}} \norm{V - V'}_\infty \}.
\end{aligned}
\end{equation}
Since $(\reals^S, \norm{\cdot}_\infty)$ is a complete metric space, $H(\reals^{S})$ is a complete metric space with respect to $d_H$. 

\begin{lem}\citep[Thm 3.3]{henrikson1999completeness}\label{lem:hausdorffComplete}
If $\mc{X}$ is a complete metric space, then its induced Hausdorff metric space $(H(\mc{X}), d_H)$ is a complete metric space.
\end{lem}
On the metric space $H(\reals^{S})$, we define a \emph{set-based Bellman operator}.
\begin{defn}[Set-based Bellman Operator]\label{def:setBellman}
For a family of MDP problems, $([S], [A], P, \mc{C}, \gamma)$, where $\mc{C} \subseteq \reals^{S\times A}$ is a compact set, its associated set-based Bellman operator is given by
\[F_{\mc{C}}(\mc{V}) = \cl\bigcup_{(C, V) \in\mc{C}\times \mc{V}} f_C(V), \quad \forall \ \mc{V} \in H(\reals^S)\]
where $\cl$ is the closure operator.
\end{defn}
As we take the union of uncountably many bounded sets, the resulting set may not be bounded, and therefore it is not immediately obvious that $F_{\mc{C}}(\mc{V})$ maps into the metric space $H(\reals^S)$. We show this is true in Proposition~\ref{prop:BellmanCompact}.
\begin{prop}\label{prop:BellmanCompact}
If $\mc{C}$ is compact, then $F_{\mc{C}}(\mc{V})\in H(\reals^{S})$, $\forall \ \mc{V}\in H(\reals^{S})$.
\end{prop}
\begin{pf}
For a non-empty set $\mc{A}$ of some finite dimensional real vector space, let us define its diameter to be denoted as $\diam{\mc{A}}=\sup_{x,y\in \mc{A}} \norm{x-y}_\infty$. The diameter of any compact set in a metric space is bounded.

We take any non-empty compact set $\mc{V}\in H(\reals^{S})$. As $F_{\mc{C}}(\mc{V})\subset \reals^S$, it suffices to prove that $F_{\mc{C}}(\mc{V})$ is closed and bounded. The closedness is guaranteed by the closure operator. A subset of a metric space is bounded iff its closure is bounded. Hence, to prove the boundedness, it suffices to prove that $\diam{\cup_{(C,V) \in\mc{C}\times\mc{V}} f_C(V)}<+\infty$. 
Consider any two cost-value function pairs, $(C,V), (C',V') \in \mc{C} \times \mc{V}$, they must satisfy
\[f_{C}(V) - f_{C'}(V') = \Big(f_{C}(V) - f_{C'}(V)\Big) + \Big(f_{C'}(V) - f_{C'}(V')\Big), \]
where the norm of the second term $\norm{f_{C'}(V) - f_{C'}(V')}_\infty$ must be upper bounded by $\norm{V - V'}_\infty$ due to contraction properties of $f_{C'}$. 
To bound the first term, we note that for any two vectors $a, b \in \reals^{S}$, \(\norm{a- b}_\infty = \max\{\max(a - b), \max(b - a)\}\) and let $\pi$ to be the optimal policy of $f_{C}(V)$,
\begin{align*}
 & \max(f_{C'}(V) - f_C(V)) \\
 \leq & \max(\nu'(\pi) + \gamma M_{\pi}P^TV - \nu(\pi)-\gamma M_{\pi}P^TV) \\ 
 \leq & \max(\nu'(\pi) - \nu(\pi)) \\
 \leq& \sum_{i \in [S]}\norm{e^T_i}_{\infty} \norm{M_{\pi}}_{\infty}\norm{ \ones_S\otimes I_A}_{\infty} \norm{(C' - C)^T}_\infty\norm{e_i}^2_{\infty}.
\end{align*}
Since $\norm{\ones_S\otimes I_A}_\infty = \norm{e_i}_\infty = \norm{e^T_i}_\infty = \norm{M_{\pi}}_\infty = 1$ for any $\pi \in \Pi$, $ \max(f_{C'}(V) - f_C(V)) \leq  S \diam{\mc{C}^T}$. The result $ \max(f_C(V) - f_{C'}(V) ) \leq  S \diam{\mc{C}^T}$  can be similarly derived. 
Therefore, $\norm{f_C(V) - f_{C'}(V')}_\infty < S\diam{\mc{C}^T} + \diam{\mc{V}}$. Since it holds for all $(C,V), (C', V') \in \mc{C} \times \mc{V}$ then $\diam{\cup_{(C,V) \in\mc{C}\times\mc{V}} f_C(V)}\leq S\diam{\mc{C}^T}+\diam{\mc{V}}<+\infty$ as $\mc{C}^T$ and $\mc{V}$ are bounded. \qed
\end{pf}
Proposition~\ref{prop:BellmanCompact} shows that $F_{\mc{C}}$ is an operator from $H(\reals^S)$ to $H(\reals^S)$. Having established its range space, we can draw many parallels between $F_{\mc{C}}$ and $f_C$. Similar to the existence of a unique fixed point $V^\star$ for $f_C$, we consider whether a \emph{fixed point set} of $F_{\mc{C}}$ which satisfies $F_{\mc{C}}(\mc{V}^\star) = \mc{V}^\star$ exists, and if it is unique. To take the comparison further, since $V^\star$ is the optimal value function for an MDP problem defined by $([S],[A],P,C,\gamma)$, how does $\mc{V}^\star$ relate to the \emph{family} of optimal solutions that corresponds to the MDP family $([S],[A],P,\mc{C}, \gamma)$? 

To prove the unique existence of $\mc{V}^\star$, we utilize the Banach fixed point theorem~\citep{puterman2014markov}, which states that a unique fixed point must exist for all contraction operators on complete metric spaces. First, we show that $F_{\mc{C}}$ is a contraction as defined in Definition~\ref{def:contractOp} on the complete metric space $(H(\reals^S), d_H)$.

\begin{prop}\label{prop:setBellman}
For any $\mc{V} \in H(\reals^S)$ and $\mc{C} \subset \reals^{S\times A}$ closed and bounded, $F_{\mc{C}}$ is a contraction operator under the Hausdorff distance. 
\end{prop}

\begin{pf}
Consider $\mc{V}$, $\bar{\mc{V}} \in H(\reals^S)$, to see that $F_{\mc{C}}$ is a contraction, we need to show 
\begin{align}
   \sup_{V \in F_{\mc{C}}(\mc{V})}\inf_{\bar{V} \in F_{\mc{C}}(\bar{\mc{V}})} \norm{V - \bar{V}}_\infty & < d_H(\mc{V}, \bar{\mc{V}}) \\
   \sup_{V \in F_{\mc{C}}(\bar{\mc{V}})} \inf_{\bar{V} \in F_{\mc{C}}(\mc{V})} \norm{V - \bar{V}}_\infty & < d_H(\mc{V}, \bar{\mc{V}}) 
\end{align}
First we note that taking $\sup$ ($\inf$) of a continuous function over a set $\mc{A}$ is equivalent to taking the $\sup$ ($\inf$) over the closure of $\mc{A}$. Let
$G_{\mc{C}}(\mc{V}) = \underset{(C,V)\in\mc{C}\times\mc{V}}{\cup} f_C(V)$ and $\cl G_{\mc{C}}(\mc{V}) =F_{\mc{C}}(\mc{V}) $, then due to continuity of norms~\citep[Thm 4.16]{rudin1964principles}, 
 \[\sup_{V\in F_{\mc{C}}(\mc{V})} \inf_{\bar{V} \in F_{\mc{C}}(\bar{\mc{V}})} \norm{V - \bar{V}} = \sup_{V\in G_{\mc{C}}(\mc{V})} \inf_{\bar{V}\in G_{\mc{C}}(\bar{\mc{V}})} \norm{V - \bar{V}}. \]

Therefore, it suffices to prove
\[{\sup_{f_C(V) \in G_{\mc{C}}(\mc{V})} \inf_{f_{\bar{C}}(\bar{V}) \in G_{\mc{C}}(\bar{\mc{V}})} \norm{f_C(V) - f_{\bar{C}}(\bar{V})}_\infty} < d_H(\mc{V}, \bar{\mc{V}}),\]
\[{\sup_{f_{\bar{C}}(\bar{V}) \in G_{\mc{C}}(\bar{\mc{V}})} \inf_{f_C(V) \in G_{\mc{C}}(\mc{V})} \norm{f_C(V) - f_{\bar{C}}(\bar{V})}_\infty } < d_H(\mc{V}, \bar{\mc{V}}).\]
For any  $V \in \mc{V}$,  $C\in\mc{C}$,
\begin{subequations}
\begin{align}
\label{prf:contraction0}
& \underset{\substack{(\bar{C}, \bar{V}) \in  \mc{C} \times\bar{\mc{V}}}}{\inf} \norm{f_C(V) - f_{\bar{C}}(\bar{V})}_\infty\\
\label{prf:contraction1}
= & \underset{\substack{(\bar{C}, \bar{V}) \in  \mc{C} \times\bar{\mc{V}}}}{\inf} \|C(\pi)+ \gamma M_{\pi}P^TV - (\bar{C}(\bar{\pi}) + \gamma M_{\bar{\pi}}P^T\bar{V})\|_\infty,\\ 
\label{prf:contraction2}
\leq & \underset{\substack{\bar{V} \in  \bar{\mc{V}}}}{\inf}  \|C(\bar{\pi}) + \gamma M_{\bar{\pi}}P^TV - (C(\bar{\pi}) + \gamma M_{\bar{\pi}}P^T\bar{V})\|_\infty,\\
\label{prf:contraction4}
\leq &  \underset{\substack{\bar{V} \in  \bar{\mc{V}}}}{\inf} \norm{\gamma M_{\bar{\pi}}P^T(V - \bar{V})}_\infty \leq  \gamma \underset{\substack{\bar{V} \in  \bar{\mc{V}}}}{\inf} \norm{V - \bar{V}}_\infty,
\end{align}
\end{subequations}
where $\pi$ corresponds to the optimal policy for the MDP $([S], [A], P, C, \gamma)$ and $\bar{\pi}$ corresponds to the optimal policy for the MDP $([S], [A], P, \bar{C}, \gamma)$ in~\eqref{prf:contraction1}. In~\eqref{prf:contraction2} we replaced $M_{\pi}$ by $M_{\bar{\pi}}$ by noting that $\pi$ is optimal, therefore $\bar{\pi}$ must result in a larger value function (similar to the proof of Prop.~\ref{prop:BellmanCompact}). In~\eqref{prf:contraction4} we note that the infimum over set $\mc{C}$ must be upper bounded by when $\bar{C} = C \in \mc{C}$, and used the fact that $\norm{M_{\bar{\pi}}P^T}_\infty \leq 1$.

Taking the $\sup$ over $G_{\mc{C}}(\mc{V})$ and $G_{\mc{C}}(\bar{\mc{V}})$, \[\sup_{V \in G_{\mc{C}}(\mc{V})}\inf_{\bar{V} \in G_{\mc{C}}(\bar{\mc{V}})} \norm{V - \bar{V}}_\infty \leq \gamma \sup_{V \in \mc{V}}\inf_{\bar{V} \in \bar{\mc{V}}} \norm{V - \bar{V}}_\infty,\] \[\sup_{\bar{V} \in G_{\mc{C}}(\bar{\mc{V}})}\inf_{V \in G_{\mc{C}}(\mc{V})} \norm{V - \bar{V}}_\infty \leq \gamma \sup_{V \in \mc{V}}\inf_{\bar{V} \in \bar{\mc{V}}} \norm{V - \bar{V}}_\infty.\] 
Therefore
$d_H(F_{\mc{C}}(\mc{V}), F_{\mc{C}}(\bar{\mc{V}})) \leq \gamma d_H(\mc{V}, \bar{\mc{V}})$. 
Since $\gamma \in (0,1)$, $F_{\mc{C}}$ is a contraction operator on $H(\reals^S)$.\qed
\end{pf}


The contraction property of $F_{\mc{C}}$ implies that repeated application of the operator to any $\mc{V}^0 \in H(\reals^S)$ will result in closer and closer sets in the Hausdorff sense of distance to a fixed point set. It is then natural to consider if there is a unique set which all $F^k_{\mc{C}}(\mc{V}^0)$ converges to. 
 
\begin{thm}\label{thm:uniqueSetFixedPont}
There exists a unique fixed point $\mc{V}^\star$ to the set-based Bellman operator $F_{\mc{C}}$ as defined in Definition~\ref{def:setBellman}, such that $F_{\mc{C}}(\mc{V}^\star) = \mc{V}^\star$, and $\mc{V}^\star$ is a closed and bounded subset of $\reals^S$. Furthermore, for any iteration starting from arbitrary $\mc{V}^0 \in H(\reals^S)$, 
$\mc{V}^{k+1} = F_{\mc{C}}(\mc{V}^k)$, the sequence converges in the Hausdorff sense i.e. $\lim_{k\to\infty} d_H(F_{\mc{C}}(\mc{V}^k),\mc{V}^\star)=0$.
\end{thm}
\begin{pf}
As shown in Proposition~\ref{prop:setBellman}, $F_{\mc{C}}$ is a contraction operator. From the Banach fixed point theorem~\citep[Thm 6.2.3]{puterman2014markov}, there exists a unique fixed point $\mc{V}^\star$, and any arbitrary $\mc{V}^0 \in H(\reals^S)$ will generate a sequence $\{F_{\mc{C}}(\mc{V}^k)\}_{k\in\naturals}$ that converges to the fixed point. \qed
\end{pf}

The fixed point $V^\star$ of Bellman operator $f_{C}$ on metric space $\reals^S$ corresponds to the optimal value function of the MDP associated with cost parameter $C$. Because there is no direct association of an MDP problem to the set-based Bellman operator $F_{\mc{C}}$, we cannot claim the same for $\mc{V}^\star$. However, $\mc{V}^\star$ does have many interesting properties on $H(\reals^S)$, in parallel to operator $f_C$ on $\reals^S$, especially in terms of the value iteration method~\eqref{eqn:VI}.  
Suppose that instead of a fixed cost parameter, we have that at each iteration $k$, a $C^k$ that is random chosen from a compact set of costs, $C^k \in \mc{C}$, then it is interesting to ask if $\mc{V}^\star$ contains all the limit points of $\lim_{k} f_{C^k}(V^k)$. Indeed, we can infer from Theorem~\ref{thm:uniqueSetFixedPont} that the sequence $\{V_k\}$ converges to $\mc{V}^\star$ under the Hausdorff metric. Furthermore, even when $V^k$ itself does not converge, it must converge to the set $\mc{V}^\star$ under the Hausdorff metric--- i.e. $\lim_{k \rightarrow 0} \inf_{V \in \mc{V}^\star}\norm{V^k - V}_\infty = 0$. 
\section{Conclusion}
We summarize our results on set-based Bellman operator: for a compact cost function set $\mc{C}$, $F_{\mc{C}}$ converges to to a unique compact set $\mc{V}^\star$ which contains all the fixed points of $f_C$ for all fixed $C \in \mc{C}$. Furthermore, $\mc{V}^\star$ also contains the limit points of $f_{C^k}(V^k)$ for any $\{C^k\}_{k\in\naturals} \subseteq \mc{C}$, $V^0 \in \reals^{S}$, given that $\lim_{k\rightarrow \infty}V^k$ converges. Even if the limit does not exist, $V^k$ must asymptotically converge to $\mc{V}^\star$ in the Hausdorff sense. Future work includes extending the uncertainty analysis to consider uncertainty in the transition kernel to fully capture learning in a general stochastic game. 


\bibliography{reference}

\begin{thebibliography}{16}
\providecommand{\natexlab}[1]{#1}
\providecommand{\url}[1]{\texttt{#1}}
\providecommand{\urlprefix}{URL }
\expandafter\ifx\csname urlstyle\endcsname\relax
  \providecommand{\doi}[1]{doi:\discretionary{}{}{}#1}\else
  \providecommand{\doi}{doi:\discretionary{}{}{}\begingroup
  \urlstyle{rm}\Url}\fi

\bibitem[{Abbad and Filar(1992)}]{abbad1992perturbation}
Abbad, M. and Filar, J.A. (1992).
\newblock Perturbation and stability theory for markov control problems.
\newblock \emph{IEEE Trans. Autom. Control}.

\bibitem[{A{\c{c}}ikme{\c{s}}e and Bayard(2012)}]{accikmecse2012markov}
A{\c{c}}ikme{\c{s}}e, B. and Bayard, D.S. (2012).
\newblock A markov chain approach to probabilistic swarm guidance.
\newblock In \emph{Amer. Control Conf.}, 6300--6307. IEEE.

\bibitem[{Altman and Gaitsgory(1993)}]{altman1993stability}
Altman, E. and Gaitsgory, V.A. (1993).
\newblock Stability and singular perturbations in constrained markov decision
  problems.
\newblock \emph{IEEE Trans. Autom. Control}, 38(6), 971--975.

\bibitem[{Ashok et~al.(2017)Ashok, Chatterjee, Daca, K{\v{r}}et{\'\i}nsk{\`y},
  and Meggendorfer}]{ashok2017value}
Ashok, P., Chatterjee, K., Daca, P., K{\v{r}}et{\'\i}nsk{\`y}, J., and
  Meggendorfer, T. (2017).
\newblock Value iteration for long-run average reward in markov decision
  processes.
\newblock In \emph{Int. Conf. Comput. Aided Verification}, 201--221. Springer.

\bibitem[{Bielecki and Filar(1991)}]{bielecki1991singularly}
Bielecki, T.R. and Filar, J.A. (1991).
\newblock Singularly perturbed markov control problem: Limiting average cost.
\newblock \emph{Ann. Op. Res.}, 28(1), 153--168.

\bibitem[{Delage and Mannor(2010)}]{delage2010percentile}
Delage, E. and Mannor, S. (2010).
\newblock Percentile optimization for markov decision processes with parameter
  uncertainty.
\newblock \emph{Op. Res.}, 58(1), 203--213.

\bibitem[{Demir et~al.(2015)Demir, Eren, and
  A{\c{c}}ikme{\c{s}}e}]{demir2015decentralized}
Demir, N., Eren, U., and A{\c{c}}ikme{\c{s}}e, B. (2015).
\newblock Decentralized probabilistic density control of autonomous swarms with
  safety.
\newblock \emph{Auton. Robots}, 39(4), 537 --554.

\bibitem[{Eisentraut et~al.(2019)Eisentraut, K{\v{r}}et{\'\i}nsk{\`y}, and
  Rotar}]{eisentraut2019stopping}
Eisentraut, J., K{\v{r}}et{\'\i}nsk{\`y}, J., and Rotar, A. (2019).
\newblock Stopping criteria for value and strategy iteration on concurrent
  stochastic reachability games.
\newblock \emph{arXiv preprint arXiv:1909.08348}.

\bibitem[{El~Chamie et~al.(2018)El~Chamie, Yu, A{\c{c}}{\i}kme{\c{s}}e, and
  Ono}]{el2018controlled}
El~Chamie, M., Yu, Y., A{\c{c}}{\i}kme{\c{s}}e, B., and Ono, M. (2018).
\newblock Controlled markov processes with safety state constraints.
\newblock \emph{IEEE Trans. Autom. Control}, 64(3), 1003--1018.

\bibitem[{Filar and Vrieze(2012)}]{filar2012competitive}
Filar, J. and Vrieze, K. (2012).
\newblock \emph{Competitive Markov decision processes}.
\newblock Springer Science \& Business Media.

\bibitem[{Givan et~al.(2000)Givan, Leach, and Dean}]{givan2000bounded}
Givan, R., Leach, S., and Dean, T. (2000).
\newblock Bounded-parameter markov decision processes.
\newblock \emph{Artif. Intell.}, 122(1-2), 71--109.

\bibitem[{Haddad and Monmege(2018)}]{haddad2018interval}
Haddad, S. and Monmege, B. (2018).
\newblock Interval iteration algorithm for mdps and imdps.
\newblock \emph{Theor. Comput. Sci.}, 735, 111--131.

\bibitem[{Henrikson(1999)}]{henrikson1999completeness}
Henrikson, J. (1999).
\newblock Completeness and total boundedness of the hausdorff metric.
\newblock In \emph{MIT Undergraduate J. Math.} Citeseer.

\bibitem[{Li et~al.(2019)Li, Yu, Calderone, Ratliff, and
  A{\c{c}}ikme{\c{s}}e}]{li2019tolling}
Li, S.H.Q., Yu, Y., Calderone, D., Ratliff, L., and A{\c{c}}ikme{\c{s}}e, B.
  (2019).
\newblock Tolling for constraint satisfaction in markov decision process
  congestion games.
\newblock In \emph{Amer. Control Conf.}, 1238--1243. IEEE.

\bibitem[{Puterman(2014)}]{puterman2014markov}
Puterman, M.L. (2014).
\newblock \emph{Markov Decision Processes: Discrete Stochastic Dynamic
  Programming}.
\newblock John Wiley \& Sons.

\bibitem[{Rudin et~al.(1964)}]{rudin1964principles}
Rudin, W. et~al. (1964).
\newblock \emph{Principles of mathematical analysis}, volume~3.
\newblock McGraw-hill New York.

\end{thebibliography}
                                                   







\end{document}